\documentclass[12pt]{article}
\usepackage{amssymb}
\usepackage{amsfonts}
\usepackage{amsmath}

\input{tcilatex}
\begin{document}

\begin{center}
\textbf{Invariant measures of Markov operators associated to iterated
function systems consisting of }${\small \varphi }$\textbf{-max-contractions
with probabilities}

\bigskip

by \textit{Flavian GEORGESCU, Radu MICULESCU }and\textit{\ Alexandru MIHAIL}

\bigskip
\end{center}

\textbf{Abstract}. {\small We prove that the Markov operator associated to
an iterated function system consisting of }${\small \varphi }${\small %
-max-contractions with probabilities has a unique invariant measure whose
support is the attractor of the system.}

\bigskip

\textbf{2010 Mathematics Subject Classification}: {\small 28A80, 37C70, 54H20%
}

\textbf{Key words and phrases}: $\varphi ${\small -max-contraction,
comparison function, iterated function system with probabilities, Markov
operator, fixed point, invariant measure}

\bigskip

\textbf{1. Introduction}

\bigskip

Iterated function systems with probabilities, which can be viewed as
particular cases of random systems with complete connections (see [8], [21],
and [24]), are well known for their applications in image compression or in
learning theory (see [1], [2], [4], [5] and the references therein).

The problem of asymptotic stability of iterated function systems consisting
of contractions with probabilities has collected a lot of attention in the
last two decades (see [9], [10], [11], [26], [28] and the references
therein).

The uniqueness of invariant probability measures for place-dependent random
iterations is treated in [19] and [25].

The problem of the existence and uniqueness of invariant measures for Markov
type operators associated to iterated function systems with probabilities
(which was initiated by J. Hutchinson [7]) was also studied, in more general
settings, in [3], [12], [14], [15], [17], [18], [23] and [27].

Since in one of our previous works we introduced a new kind of iterated
function systems, namely those consisting of $\varphi $-max-contractions,
and we prove the existence and uniqueness of their attractor (see [6]),
along the lines of research previously mentioned, the next step, -which is
accomplished in the present paper- is to study the Markov operators
associated to such systems with probabilities. We prove that each such
operator has a unique invariant measure whose support is the attractor of
the system. Let us point out that the invariant measure is obtained via the
Riesz representation theorem from a positive linear functional \ which is
generated by the dual operator of the Markov operator.

\bigskip

\textbf{2.} \textbf{Preliminaries}

\bigskip

\textbf{Notations and terminology}

\bigskip

Given the sets $A$ and $B$, by $B^{A}$ we mean the set of functions from $A$
to $B$.

Given a set $X$, a function $f:X\rightarrow X$ and $n\in \mathbb{N}$, by $%
f^{[n]}$ we mean $\underset{n\text{ times}}{f\circ f\circ ...\circ f}$.

Given a metric space $(X,d)$, by:

- $diam(A)$ we mean the diameter of the subset $A$ of $X$

- $P_{cp}(X)$ we mean the set of non-empty compact subsets of $X$

- $\mathcal{C}(X)$ we mean the set of continuous functions $f:X\rightarrow 
\mathbb{R}$

- $\mathcal{B}(X)$ we mean the $\sigma $-algebra of Borel subsets of $X$

- the support of a finite positive borelian measure $\mu $ on $X$ (denoted
by supp $\mu $) we mean the smallest closed subset of $X$ on which $\mu $ is
concentrated; so 
\begin{equation*}
\text{supp }\mu =\underset{F=\overline{F}\subseteq X\text{, }\mu (F)=\mu (X)}%
{\cap }F
\end{equation*}%
\smallskip

- $\mathcal{M}(X)$ we mean the space of borelian normalized and positive
measures on $X$ with compact support

- $Lip_{1}(X,\mathbb{R})$ we mean the set of functions $f:X\rightarrow 
\mathbb{R}$ having the property that $lip(f)\overset{def}{=}\underset{x,y\in
X,x\neq y}{\sup }\frac{\left\vert f(x)-f(y)\right\vert }{d(x,y)}\leq 1$

- the Hausdorff-Pompeiu metric we mean\textit{\ }$H:P_{cp}(X)\times
P_{cp}(X)\rightarrow \lbrack 0,+\infty )$\textit{\ }given by

\begin{equation*}
H(A,B)=\max \{\underset{x\in A}{\sup }(\underset{y\in B}{\inf }d(x,y)),%
\underset{x\in B}{\sup }(\underset{y\in A}{\inf }d(x,y))\}
\end{equation*}%
for all $A,B\in P_{cp}(X)$

- a Picard operator we mean a function $f:X\rightarrow X$\textit{\ }having
the property that there exists a unique fixed point $\alpha $ of $f$\ and
the sequence $(f^{[n]}(x))_{n\in \mathbb{N}}$\ is convergent to $\alpha $
for every $x\in X$.

\newpage

\textbf{The Hutchinson distance}

\bigskip

\textbf{Definition 2.1.} \textit{Given a complete metric space }$(X,d)$%
\textit{, the function} $d_{H}:\mathcal{M}(X)\times \mathcal{M}%
(X)\rightarrow \lbrack 0,\infty )$ \textit{described by} $d_{H}(\mu ,\nu )=%
\underset{f\in Lip_{1}(X,\mathbb{R})}{\sup }\left\vert \underset{X}{\dint }%
fd\mu -\underset{X}{\dint }fd\nu \right\vert $\textit{\ for every} $\mu ,\nu
\in \mathcal{M}(X)$\textit{, turns out to be a distance which is called the
Hutchinson distance.}

\bigskip

\textbf{Remark 2.2} (see page 46 from [22])\textbf{.} \textit{Given a
compact metric space} $(X,d)$, $\mu \in \mathcal{M}(X)$ \textit{and a
sequence} $(\mu _{n})_{n\in \mathbb{N}}$ \textit{of elements from} $\mathcal{%
M}(X)$\textit{, the following statements are equivalent:}

\textit{a) the sequence }$(\mu _{n})_{n\in \mathbb{N}}$\textit{\ converges
to }$\mu $\textit{\ in the weak topology i.e. }$\underset{X}{\dint }gd\mu =%
\underset{n\rightarrow \infty }{\lim }\underset{X}{\dint }gd\mu _{n}$\textit{%
\ for every }$g\in \mathcal{C}(X)$;

\textit{b)} $\underset{n\rightarrow \infty }{\lim }d_{H}(\mu _{n},\mu )=0$%
\textit{.}

\bigskip

\textbf{Comparison functions}

\bigskip

\textbf{Definition 2.3.} \textit{A function }$\varphi :[0,\infty
)\rightarrow \lbrack 0,\infty )$\textit{\ is called a comparison function
provided that it satisfies the following properties:}

\textit{i) }$\varphi $\textit{\ is increasing;}

\textit{ii) }$\underset{n\rightarrow \infty }{\lim }\varphi ^{\lbrack
n]}(x)=0$\textit{\ for every }$x\in \lbrack 0,\infty )$\textit{.}

\bigskip

\textbf{Remark 2.4.} \textit{For each comparison function the following two
properties are valid:}

\textit{a) }$\varphi (0)=0$\textit{;}

\textit{b) }$\varphi (x)<x$ \textit{for every }$x\in (0,\infty )$\textit{.}

\bigskip

\textbf{The shift space}

\bigskip

Given a nonempty set $I$, we denote the set $I^{\mathbb{N}^{\ast }}$ by $%
\Lambda (I)$. Thus $\Lambda (I)$ is the set of infinite words with letters
from the alphabet $I$ and a standard element $\omega $ of $\Lambda (I)$ can
be presented as $\omega =\omega _{1}\omega _{2}...\omega _{n}\omega
_{n+1}... $ .

Given a nonempty set $I$, we denote the set $I^{\{1,2,...,n\}}$ by $\Lambda
_{n}(I)$. Thus $\Lambda _{n}(I)$ is the set of words of length $n$ with
letters from the alphabet $I$ and a standard element $\omega $ of $\Lambda
_{n}(I)$ can be presented as $\omega =\omega _{1}\omega _{2}...\omega _{n}$.
By $\Lambda _{0}(I)$ we mean the set having only one element, namely the
empty word denoted by $\lambda $.

For $n\in \mathbb{N}^{\ast }$, we shall use the following notation: $V_{n}(I)%
\overset{not}{=}\underset{k\in \{0,1,2,...,n-1\}}{\cup }\Lambda _{k}(I)$.

Given a nonempty set $I$, $m,n\in \mathbb{N}$ and two words $\omega =\omega
_{1}\omega _{2}...\omega _{n}\in \Lambda _{n}(I)$\ and $\theta =\theta
_{1}\theta _{2}...\theta _{m}\in \Lambda _{m}(I)$ or $\theta =\theta
_{1}\theta _{2}...\theta _{m}\theta _{m+1}...\in \Lambda (I)$, by $\omega
\theta $ we mean the concatenation of the words $\omega $ and $\theta $, i.e.%
$\ \omega \theta =\omega _{1}\omega _{2}...\omega _{n}\theta _{1}\theta
_{2}...\theta _{m}$ and respectively $\omega \theta =\omega _{1}\omega
_{2}...\omega _{n}\theta _{1}\theta _{2}...\theta _{m}\theta _{m+1}...$ .

For a family of functions $(f_{i})_{i\in I}$, where $f_{i}:X\rightarrow X$,
and $\omega =\omega _{1}\omega _{2}...\omega _{n}\in \Lambda _{n}(I)$, we
shall use the following notation: $f_{\omega }\overset{not}{=}f_{\omega
_{1}}\circ f_{\omega _{2}}\circ ...\circ f_{\omega _{n}}$.

For a function $f:X\rightarrow X$, by $f_{\lambda }$ we mean $Id_{X}$.

\bigskip

\textbf{A result concerning a sequence of compact subsets of a metric space}

\bigskip

\textbf{Proposition 2.5 }(see Proposition 2.8 from [16])\textbf{.} \textit{%
Let }$(X,d)$\textit{\ be a complete metric space, }$(Y_{n})_{n\in \mathbb{N}%
}\subseteq P_{cp}(X)$ \textit{and }$Y\in P_{cp}(X)$\textit{\ such that }$%
\underset{n\rightarrow \infty }{\lim }H(Y_{n},Y)=0$\textit{. Then }$Y\cup (%
\overset{\infty }{\underset{n=0}{\cup }}Y_{n})\in P_{cp}(X)$\textit{.}

\bigskip

\textbf{3.} \textbf{The Markov operator associated to an}\textit{\ }\textbf{%
iterated function system consisting of }$\varphi $-\textbf{max-contractions
with probabilities }

\bigskip

\textbf{Definition 3.1.\ }\textit{An iterated function system consisting of }%
$\varphi $-$\max $\textit{-contractions} \textit{(}$\varphi $-$\max $-%
\textit{IFS for short) is described by:}

-\textit{\ a complete metric space }$(X,d)$

-\textit{\ a finite family of continuous functions }$(f_{i})_{i\in I}$%
\textit{, where }$f_{i}:X\rightarrow X$,\textit{\ having the property that
there exist a comparison function }$\varphi :[0,\infty )\rightarrow \lbrack
0,\infty )$\textit{\ and }$p\in \mathbb{N}^{\ast }$ \textit{such that} $%
\underset{\omega \in \Lambda _{p}(I)}{\max }d(f_{\omega }(x),f_{\omega
}(y))\leq \varphi (\underset{\omega \in V_{p}(I)}{\max }d(f_{\omega
}(x),f_{\omega }(y)))$ \textit{for every }$x,y\in X$\textit{.}

\medskip

\textit{We denote such a system by }$\mathcal{S}=((X,d),(f_{i})_{i\in I})$%
\textit{.}

\medskip

The \textit{fractal operator} $F_{\mathcal{S}}:P_{cp}(X)\rightarrow
P_{cp}(X) $, associated to the $\varphi $-$\max $-IFS $\mathcal{S}$, is
given by $F_{\mathcal{S}}(K)=\underset{i\in I}{\cup }f_{i}(K)$ for every $%
K\in P_{cp}(X)$.

\medskip

We say that the $\varphi $-$\max $-IFS $\mathcal{S}$ has attractor if $F_{%
\mathcal{S}}$ is a Picard operator (with respect to the Hausdorff-Pompeiu
metric) and the fixed point of $F_{\mathcal{S}}$ is called the \textit{%
attractor} of the system $\mathcal{S}$ and it is denoted by $A_{\mathcal{S}}$%
.

\bigskip

\textbf{Theorem 3.2 }(see Theorem 3.2 from [6])\textbf{.} \textit{Each} $%
\varphi $-$\max $\textit{-IFS has attractor.}

\bigskip

\textbf{Definition 3.3.\ }\textit{An iterated function system consisting of }%
$\varphi $-$\max $\textit{-contractions with probabilities} \textit{(}$%
\varphi $-$\max $-\textit{IFSp for short) is described by:}

\textit{- a }$\varphi $-$\max $\textit{-IFS }$\mathcal{S}=((X,d),(f_{i})_{i%
\in \{1,2,...,m\}})$

\textit{- a system of probabilities }$(p_{i})_{i\in \{1,2,...,m\}}$\textit{,
i.e. }$p_{i}\in (0,1)$\textit{\ for every }$i\in \{1,2,...,m\}$\textit{\ and 
}$p_{1}+p_{2}+...+p_{m}=1$\textit{.}

\medskip

\textit{We denote such a system by }$\mathcal{S}=((X,d),(f_{i})_{i\in
\{1,2,...,m\}},(p_{i})_{i\in \{1,2,...,m\}})$\textit{.}

\medskip

We associate to such a system \textit{the Markov operator} $M_{\mathcal{S}}:%
\mathcal{M}(X)\rightarrow \mathcal{M}(X)$ given by $M_{\mathcal{S}}(\mu
)=p_{1}\mu \circ f_{1}^{-1}+...+p_{m}\mu \circ f_{m}^{-1}$, i.e. $M_{%
\mathcal{S}}(\mu )(B)=p_{1}\mu (f_{1}^{-1}(B))+...+p_{m}\mu (f_{m}^{-1}(B))$
for every $B\in \mathcal{B}(X)$ and every $\mu \in \mathcal{M}(X)$. A fixed
point of $M_{\mathcal{S}}$ is called \textit{invariant measure}.

\bigskip

\textbf{Lemma 3.4} (see Lemma 3 from [9])\textbf{. }\textit{The equality}
supp $M_{\mathcal{S}}^{[n]}(\mu )=F_{\mathcal{S}}^{[n]}($supp $\mu )$\ 
\textit{is valid} \textit{for every }$\varphi $-$\max $-\textit{IFSp }$%
\mathcal{S}$\textit{, every} $\mu \in \mathcal{M}(X)$ \textit{and every} $%
n\in \mathbb{N}$\textit{.}

\bigskip

\textbf{Remark 3.5.}\textit{\ The operator }$M_{\mathcal{S}}$ \textit{is
well defined, for every} $\varphi $-$\max $-\textit{IFSp }$\mathcal{S}$.

\bigskip

\textbf{Remark 3.6.} \textit{Given a} $\varphi $-$\max $-\textit{IFSp }$%
\mathcal{S}=((X,d),(f_{i})_{i\in \{1,...,m\}},(p_{i})_{i\in \{1,...,m\}})$%
\textit{, the Markov operator} $M_{\mathcal{S}}$ \textit{is a Feller operator%
} \textit{since} $\underset{X}{\dint }gdM_{\mathcal{S}}(\mu )=p_{1}\underset{%
X}{\dint }g\circ f_{1}d\mu $\linebreak $+...+p_{m}\underset{X}{\dint }g\circ
f_{m}d\mu $ \textit{for every continuous and bounded function} $%
g:X\rightarrow \mathbb{R}$\textit{\ and every }$\mu \in \mathcal{M}(X)$%
\textit{.}

\bigskip

\textbf{4. The main result}

\bigskip

Our main result states that the Markov operator associated to an IFSp is a
Picard operator and the support of its fixed point is the attractor of the
system. At the beginning, we shall consider for the case of a system for
which the metric space is compact (see Theorem 4.9) and then we shall treat
the general case (see Theorem 4.18).

\bigskip

\textbf{A. The case of a} $\varphi $-$\max $\textbf{-IFSp for which the
metric space is compact}

\bigskip

Let us start with some:

\medskip

\textbf{Notations}. For a $\varphi $-$\max $-IFSp $\mathcal{S}%
=((X,d),(f_{i})_{i\in \{1,2,...,m\}},(p_{i})_{i\in \{1,2,...,m\}})$\textit{,}
$g:X\rightarrow \mathbb{R}$, $x,y\in X$, $n\in \mathbb{N}$, $\omega =\omega
_{1}\omega _{2}...\omega _{n}\in \Lambda _{n}(\{1,2,...,m\})$ and\textit{\ }$%
\varepsilon >0$, we shall use the following notations:

$\bullet $ $X_{x,y,n}\overset{\text{not}}{=}\{d(f_{\omega }(x),f_{\omega
}(y))\mid \omega \in \Lambda _{n}(\{1,2,...,m\})\}$

$\bullet $ $X_{\varepsilon ,n}\overset{\text{not}}{=}\{\max X_{x,y,n}\mid
x,y\in X$, $d(x,y)\leq \varepsilon \}$

$\bullet $ $a_{n}(\varepsilon )\overset{\text{not}}{=}\sup X_{\varepsilon
,n} $

$\bullet $ $\mathcal{O}_{\varepsilon }(g)\overset{\text{not}}{=}\underset{%
x,y\in X\text{, }d(x,y)\leq \varepsilon }{\sup }\left\vert
g(x)-g(y)\right\vert $

$\bullet $ $B_{\mathcal{S}}(g)\overset{\text{not}}{=}p_{1}g\circ
f_{1}+p_{2}g\circ f_{2}+...+p_{m}g\circ f_{m}$

$\bullet $ $p_{\omega }\overset{\text{not}}{=}p_{\omega _{1}}p_{\omega
_{2}}...p_{\omega _{n}}$.

\bigskip

An easy mathematical induction argument proves the following:

\bigskip

\textbf{Remark 4.1.} \textit{For every} $\varphi $-$\max $-\textit{IFSp} $%
\mathcal{S}=((X,d),(f_{i})_{i\in \{1,...,m\}},(p_{i})_{i\in \{1,...,m\}})$%
\textit{,} $g:X\rightarrow \mathbb{R}$ \textit{and }$n\in \mathbb{N}$\textit{%
, we have} $B_{\mathcal{S}}^{[n]}(g)=\underset{\omega \in \Lambda
_{n}(\{1,2,...,m\})}{\sum }p_{\omega }g\circ f_{\omega }$.

\bigskip

\textbf{Lemma 4.2.} \textit{For every} $\varphi $-$\max $-\textit{IFSp} $%
\mathcal{S}$\textit{,} $g:X\rightarrow \mathbb{R}$, $n\in \mathbb{N}$ 
\textit{and }$\varepsilon >0$,\textit{\ we have} $\mathcal{O}_{\varepsilon
}(B_{\mathcal{S}}^{[n]}(g))\leq \mathcal{O}_{a_{n}(\varepsilon )}(g)$.

\textit{Proof}. Let us suppose that $\mathcal{S}=((X,d),(f_{i})_{i\in
\{1,2,...,m\}},(p_{i})_{i\in \{1,2,...,m\}})$. Then we have 
\begin{equation*}
\mathcal{O}_{\varepsilon }(B_{\mathcal{S}}^{[n]}(g))=\underset{x,y\in X\text{%
, }d(x,y)\leq \varepsilon }{\sup }\left\vert B_{\mathcal{S}}^{[n]}(g)(x)-B_{%
\mathcal{S}}^{[n]}(g)(y)\right\vert \overset{\text{Remark 4.1}}{=}
\end{equation*}%
\begin{equation*}
=\underset{x,y\in X\text{, }d(x,y)\leq \varepsilon }{\sup }\left\vert 
\underset{\omega \in \Lambda _{n}(\{1,2,...,m\})}{\sum }p_{\omega }((g\circ
f_{\omega })(x)-(g\circ f_{\omega })(y))\right\vert \leq 
\end{equation*}%
\begin{equation*}
\leq \underset{x,y\in X\text{, }d(x,y)\leq \varepsilon }{\sup }\underset{%
\omega \in \Lambda _{n}(\{1,2,...,m\})}{\sum }p_{\omega }\left\vert
g(f_{\omega }(x))-g(f_{\omega }(y))\right\vert \leq 
\end{equation*}%
\begin{equation*}
\leq \underset{\omega \in \Lambda _{n}(\{1,2,...,m\})}{\sum }p_{\omega }%
\mathcal{O}_{a_{n}(\varepsilon )}(g)=\mathcal{O}_{a_{n}(\varepsilon )}(g)%
\underset{\omega \in \Lambda _{n}(\{1,2,...,m\})}{\sum }p_{\omega }=\mathcal{%
O}_{a_{n}(\varepsilon )}(g)\text{. }\square 
\end{equation*}

\bigskip

\textbf{Lemma 4.3.} \textit{For every} $\varphi $-$\max $-\textit{IFSp} $%
\mathcal{S}$\textit{,} $n\in \mathbb{N}$ \textit{and }$\varepsilon >0$%
\textit{, we have} $a_{n+p}(\varepsilon )\leq \varphi (\max
\{a_{n+p-1}(\varepsilon ),a_{n+p-2}(\varepsilon ),...,a_{n}(\varepsilon )\})$%
, \textit{where the meaning of the natural number }$p$\textit{\ is the one
from Definition 3.1.}

\textit{Proof}. Let us suppose that $\mathcal{S}=((X,d),(f_{i})_{i\in
\{1,2,...,m\}},(p_{i})_{i\in \{1,2,...,m\}})$ and suppose that $\varphi $ is
the comparison function described in Definition 3.1. For the sake of
simplicity, in the framework of this proof, we denote $a_{n}(\varepsilon )$
by $a_{n}$ and $\max \{a_{n+p-1},a_{n+p-2},...,a_{n}\}$ by $M_{n}$. For all $%
x,y\in X$, $d(x,y)\leq \varepsilon $, $\theta \in V_{p}(\{1,2,...,m\})$ and $%
\omega \in \Lambda _{n}(\{1,2,...,m\})$, we have $d(f_{\theta }(f_{\omega
}(x)),f_{\theta }(f_{\omega }(y)))\leq a_{n+\left\vert \theta \right\vert
}\leq M_{n}$. Thus $\underset{\theta \in V_{p}(\{1,2,...,m\})}{\max }%
d(f_{\theta }(f_{\omega }(x)),f_{\theta }(f_{\omega }(y)))\leq M_{n}$, so%
\begin{equation}
\varphi (\underset{\theta \in V_{p}(\{1,2,...,m\})}{\max }d(f_{\theta
}(f_{\omega }(x)),f_{\theta }(f_{\omega }(y))))\leq \varphi (M_{n})\text{,} 
\tag{1}
\end{equation}%
for every $x,y\in X$, $d(x,y)\leq \varepsilon $ and $\omega \in \Lambda
_{n}(\{1,2,...,m\})$.

For each $\omega ^{^{\prime }}=\omega _{1}\omega _{2}...\omega _{p-1}\omega
_{p}\omega _{p+1}...\omega _{n+p}\in \Lambda _{n+p}(\{1,2,...,m\})$, $x,y\in
X$, $d(x,y)\leq \varepsilon $, with the notations $v\overset{not}{=}\omega
_{1}\omega _{2}...\omega _{p-1}\omega _{p}$ and $w\overset{not}{=}\omega
_{p+1}...\omega _{n+p}$, we have $d(f_{\omega }(x)),f_{\omega }(y))\overset{%
\text{Definition 3.1}}{\leq }\varphi (\underset{\theta \in
V_{p}(\{1,2,...,m\})}{\max }d(f_{\theta }(f_{w}(x)),f_{\theta }(f_{w}(y))))%
\overset{(1)}{\leq }\varphi (M_{n})$. Hence $\max X_{x,y,n+p}\leq \varphi
(M_{n})$ for every $x,y\in X$, $d(x,y)\leq \varepsilon $. Therefore, we come
to the conclusion that $\underset{x,y\in X,d(x,y)\leq \varepsilon }{\sup }%
\max X_{x,y,n+p}\leq \varphi (M_{n})$, i.e. $a_{n+p}\leq \varphi (\max
\{a_{n+p-1},a_{n+p-2},...,a_{n}\})$. $\square $

\bigskip

\textbf{Lemma 4.4.} \textit{For every} $\varphi $-$\max $-\textit{IFSp} $%
\mathcal{S}$ \textit{and }$\varepsilon >0$\textit{, we have} $\underset{%
n\rightarrow \infty }{\lim }a_{n}(\varepsilon )=0$.

\textit{Proof}. Let us suppose that $\mathcal{S}=((X,d),(f_{i})_{i\in
\{1,2,...,m\}},(p_{i})_{i\in \{1,2,...,m\}})$ and the meaning of $p$ and $%
\varphi $ is the one described in Definition 3.1. For the sake of
simplicity, in the framework of this proof, we denote $a_{n}(\varepsilon )$
by $a_{n}$ and $\max \{a_{1},a_{2},...,a_{p}\}$ by $M$.

\medskip

\textbf{Claim}. $0\leq a_{pk+j}\leq \varphi ^{\lbrack k]}(M)$ \textit{for
every} $k\in \mathbb{N}$ \textit{and every\ }$j\in \{1,2,...,p\}$.

Indeed, we have 
\begin{equation}
a_{p+1}\overset{\text{Lemma 4.3}}{\leq }\varphi (M)\overset{\text{Remark 2.4}%
}{\leq }M\text{.}  \tag{1}
\end{equation}%
We also have 
\begin{equation}
a_{p+2}\overset{\text{Lemma 4.3}}{\leq }\varphi (\max
\{a_{p+1},a_{p},...,a_{2}\})\text{.}  \tag{2}
\end{equation}%
The inequalities $a_{2}\leq M,...,a_{p}\leq M$ and $(1)$ lead to the
conclusion that $\max \{a_{p+1},a_{p},...,a_{2}\}\leq M$ and in view of $(2)$
we get $a_{p+2}\leq \varphi (M)$ and continuing the same line of reasoning
we obtain that 
\begin{equation}
a_{p+3}\leq \varphi (M)\text{, }...\text{, }a_{2p}\leq \varphi (M)\text{.} 
\tag{3}
\end{equation}%
Moreover, we have 
\begin{equation*}
a_{2p+1}\overset{\text{Lemma 4.3}}{\leq }\varphi (\max
\{a_{2p},a_{2p-1},...,a_{p+1}\})\overset{(1),(2)\&(3)}{\leq }\varphi
(\varphi (M))=\varphi ^{\lbrack 2]}(M)
\end{equation*}%
and, as above, we come to the conclusion that $a_{2p+2}\leq \varphi
^{\lbrack 2]}(M)$, $...$, $a_{3p}\leq \varphi ^{\lbrack 2]}(M)$. Now,
inductively one can prove the claim.

\medskip

Since $\underset{k\rightarrow \infty }{\lim }\varphi ^{\lbrack k]}(M)=0$
(see Definition 2.3), based on the Claim and the squeeze theorem, the proof
is done. $\square $

\bigskip

\textbf{Lemma 4.5.} \textit{For every} $\varphi $-$\max $-\textit{IFSp} $%
\mathcal{S}=((X,d),(f_{i})_{i\in \{1,...,m\}},(p_{i})_{i\in \{1,...,m\}})$, 
\textit{with }$(X,d)$\textit{\ compact,} $g:X\rightarrow \mathbb{R}$ \textit{%
continuous} \textit{and }$\varepsilon >0$,\textit{\ we have} $\underset{%
n\rightarrow \infty }{\lim }\mathcal{O}_{a_{n}(\varepsilon )}(g)$\linebreak $%
=0$.

\textit{Proof}. As $(X,d)$ is compact and $g$ is continuous, for every $%
\varepsilon ^{^{\prime }}>0$ there exists $\delta _{\varepsilon ^{^{\prime
}}}>0$ such that $\left\vert g(x)-g(y)\right\vert <\varepsilon ^{^{\prime }}$
for every $x,y\in X$ such that $d(x,y)<\delta _{\varepsilon ^{^{\prime }}}$.
According to Lemma 4.4 there exists $n_{\varepsilon ^{^{\prime }}}\in 
\mathbb{N}$ such that $a_{n}(\varepsilon )<\delta _{\varepsilon ^{^{\prime
}}}$ for every $n\in \mathbb{N}$, $n\geq n_{\varepsilon ^{^{\prime }}}$.
Consequently we get $\left\vert g(x)-g(y)\right\vert <\varepsilon ^{^{\prime
}}$ for every $x,y\in X$ such that $d(x,y)<a_{n}(\varepsilon )$ and $n\in 
\mathbb{N}$, $n\geq n_{\varepsilon ^{^{\prime }}}$, so $\underset{x,y\in X%
\text{, }d(x,y)\leq a_{n}(\varepsilon )}{\sup }\left\vert
g(x)-g(y)\right\vert \leq \varepsilon ^{^{\prime }}$, i.e. $\mathcal{O}%
_{a_{n}(\varepsilon )}(g)\leq \varepsilon ^{^{\prime }}$ for every $n\in 
\mathbb{N}$, $n\geq n_{\varepsilon ^{^{\prime }}}$. This means that $%
\underset{n\rightarrow \infty }{\lim }\mathcal{O}_{a_{n}(\varepsilon )}(g)=0$%
. $\square $

\bigskip

\textbf{Lemma 4.6.} \textit{For every} $\varphi $-$\max $-\textit{IFSp} $%
\mathcal{S}=((X,d),(f_{i})_{i\in \{1,...,m\}},(p_{i})_{i\in \{1,...,m\}})$, 
\textit{with }$(X,d)$\textit{\ compact,} $g:X\rightarrow \mathbb{R}$\textit{%
\ continuous} \textit{and }$\varepsilon >0$,\textit{\ we have} $\underset{%
n\rightarrow \infty }{\lim }\mathcal{O}_{\varepsilon }(B_{\mathcal{S}%
}^{[n]}(g))$\linebreak $=0$.

\textit{Proof}. The squeeze theorem, Lemma 4.2 and Lemma 4.5 assure us the
validity of this Lemma. $\square $

\bigskip

\textbf{Proposition 4.7.} \textit{For every} $\varphi $-$\max $-\textit{IFSp}
$\mathcal{S}=((X,d),(f_{i})_{i\in \{1,...,m\}},(p_{i})_{i\in \{1,...,m\}})$, 
\textit{with }$(X,d)$\textit{\ compact,\ and} $g:X\rightarrow \mathbb{R}$%
\textit{\ continuous, there exists a constant function} $c_{g}:X\rightarrow 
\mathbb{R}$ \textit{such that} $B_{\mathcal{S}}^{[n]}(g)\underset{%
n\rightarrow \infty }{\overset{u}{\rightarrow }}c_{g}$.

\textit{Proof}. We divide the proof into three steps.

\medskip

\textbf{Step 1}.\textit{\ The sequence} $(\underset{x\in X}{\sup }B_{%
\mathcal{S}}^{[n]}(g)(x))_{n\in \mathbb{N}}$ \textit{is decreasing and the
sequence} $(\underset{x\in X}{\inf }B_{\mathcal{S}}^{[n]}(g)(x))_{n\in 
\mathbb{N}}$ \textit{is increasing.}

We have%
\begin{equation*}
B_{\mathcal{S}}^{[n+1]}(g)(x)\overset{\text{Remark 4.1}}{=}\underset{\omega
\in \Lambda _{n+1}(\{1,2,...,m\})}{\sum }p_{\omega }(g\circ f_{\omega })(x)=
\end{equation*}%
\begin{equation*}
=\underset{i\in \{1,2,...,m\},v\in \Lambda _{n}(\{1,2,...,m\})}{\sum }%
p_{i}p_{v}(g\circ f_{v}\circ f_{i})(x)=
\end{equation*}%
\begin{equation*}
=\underset{i\in \{1,2,...,m\}}{\sum }p_{i}(\underset{v\in \Lambda
_{n}(\{1,2,...,m\})}{\sum }p_{v}(g\circ f_{v})(f_{i}(x)))\overset{\text{%
Remark 4.1}}{=}
\end{equation*}%
\begin{equation*}
=\underset{i\in \{1,2,...,m\}}{\sum }p_{i}B_{\mathcal{S}}^{[n]}(g)(f_{i}(x))%
\leq \underset{i\in \{1,2,...,m\}}{\sum }p_{i}\underset{x\in X}{\sup }B_{%
\mathcal{S}}^{[n]}(g)(x)=\underset{x\in X}{\sup }B_{\mathcal{S}}^{[n]}(g)(x)%
\text{,}
\end{equation*}%
for every $x\in X$ and every $n\in \mathbb{N}$, so $\underset{x\in X}{\sup }%
B_{\mathcal{S}}^{[n+1]}(g)(x)\leq \underset{x\in X}{\sup }B_{\mathcal{S}%
}^{[n]}(g)(x)$ for every $n\in \mathbb{N}$, i.e. the sequence $(\underset{%
x\in A_{\mathcal{S}}}{\sup }B_{\mathcal{S}}^{[n]}(g)(x))_{n\in \mathbb{N}}$
is decreasing. In a similar manner one can prove that the sequence $(%
\underset{x\in X}{\inf }B_{\mathcal{S}}^{[n]}(g)(x))_{n\in \mathbb{N}}$ is
increasing.

\medskip

\textbf{Step 2}.\textit{\ The sequences} $(\underset{x\in X}{\inf }B_{%
\mathcal{S}}^{[n]}(g)(x))_{n\in \mathbb{N}}$ \textit{and} $(\underset{x\in X}%
{\sup }B_{\mathcal{S}}^{[n]}(g)(x))_{n\in \mathbb{N}}$ \textit{are
convergent and they have the same limit (which will be denoted by }$c_{g}$%
\textit{).}

Step 1 assures us that there exist $l_{1},l_{2}\in \mathbb{R}$ such that $%
\underset{n\rightarrow \infty }{\lim }\underset{x\in X}{\inf }B_{\mathcal{S}%
}^{[n]}(g)(x)=l_{2}\leq l_{1}=\underset{n\rightarrow \infty }{\lim }\underset%
{x\in X}{\sup }B_{\mathcal{S}}^{[n]}(g)(x)$. We have $0\leq l_{1}-l_{2}\leq 
\underset{x\in X}{\sup }B_{\mathcal{S}}^{[n]}(g)(x)-\underset{x\in X}{\inf }%
B_{\mathcal{S}}^{[n]}(g)(x)=\underset{x,y\in X}{\sup }(B_{\mathcal{S}%
}^{[n]}(g)(x)-B_{\mathcal{S}}^{[n]}(g)(y))\leq \mathcal{O}_{diam(X)}B_{%
\mathcal{S}}^{[n]}(g)$ for every $n\in \mathbb{N}$ and by passing to limit
as $n\rightarrow \infty $, based on Lemma 4.6, we get $l_{1}=l_{2}\overset{%
not}{=}c_{g}$.

\medskip

\textbf{Step 3}.\textit{\ There exists a constant function} $%
c_{g}:X\rightarrow \mathbb{R}$ \textit{such that} $B_{\mathcal{S}}^{[n]}(g)%
\underset{n\rightarrow \infty }{\overset{u}{\rightarrow }}c_{g}$.

Considering the constant function $c_{g}:X\rightarrow \mathbb{R}$ given by $%
c_{g}(x)=c_{g}$ for every $x\in X$, we have $-\mathcal{O}_{diam(X)}B_{%
\mathcal{S}}^{[n]}(g)=\underset{x\in X}{\inf }B_{\mathcal{S}}^{[n]}(g)(x)-%
\underset{x\in X}{\sup }B_{\mathcal{S}}^{[n]}(g)(x)\leq \underset{x\in X}{%
\inf }B_{\mathcal{S}}^{[n]}(g)(x)-c_{g}\leq B_{\mathcal{S}%
}^{[n]}(g)(x)-c_{g}\leq \underset{x\in X}{\sup }B_{\mathcal{S}}^{[n]}(g)(x)-%
\underset{x\in X}{\inf }B_{\mathcal{S}}^{[n]}(g)(x)=\mathcal{O}_{diam(X)}B_{%
\mathcal{S}}^{[n]}(g)$, i.e. $\left\vert B_{\mathcal{S}%
}^{[n]}(g)(x)-c_{g}(x)\right\vert \leq \mathcal{O}_{diam(X)}B_{\mathcal{S}%
}^{[n]}(g)$, for every $x\in X$ and every $n\in \mathbb{N}$. The last
inequality and Lemma 4.6 assure us that $B_{\mathcal{S}}^{[n]}(g)\underset{%
n\rightarrow \infty }{\overset{u}{\rightarrow }}c_{g}$. $\square $

\bigskip

\textbf{Proposition 4.8.} \textit{For every} $\varphi $-$\max $-\textit{IFSp}
$\mathcal{S}=((X,d),(f_{i})_{i\in \{1,...,m\}},(p_{i})_{i\in \{1,...,m\}})$, 
\textit{with }$(X,d)$\textit{\ compact,} \textit{there exists a unique} 
\textit{borelian positive measure }$\mu _{\mathcal{S}}$\textit{\ on }$X$ 
\textit{such that} $c_{g}=\underset{X}{\int }gd\mu _{\mathcal{S}}$ \textit{%
for every continuous function }$g:X\rightarrow \mathbb{R}$\textit{.}

\textit{Proof}. Let us consider the function $I:\mathcal{C}(X)\rightarrow 
\mathbb{R}$ given by $I(g)=c_{g}$ for every $g\in \mathcal{C}(X)$. As,
according to Remark 4.1, $B_{\mathcal{S}}^{[n]}(g+h)=B_{\mathcal{S}%
}^{[n]}(g)+B_{\mathcal{S}}^{[n]}(h)$\textit{\ }and $B_{\mathcal{S}%
}^{[n]}(\alpha g)=\alpha B_{\mathcal{S}}^{[n]}(g)$ for every $g,h\in 
\mathcal{C}(X)$, $\alpha \in \mathbb{R}$ and every $n\in \mathbb{N}$, by
passing to limit as $n\rightarrow \infty $, we get $I(g+h)=I(g)+I(h)$ and $%
I(\alpha g)=\alpha I(g)$ for every $g,h\in \mathcal{C}(X)$ and every $\alpha
\in \mathbb{R}$. Moreover, as $B_{\mathcal{S}}^{[n]}(g)\geq 0$ for every $%
n\in \mathbb{N}$ and every $g\in \mathcal{C}(X)$ such that $g\geq 0$, by
passing to limit as $n\rightarrow \infty $, we get $I(g)\geq 0$ for every $%
g\in \mathcal{C}(X)$ such that $g\geq 0$. We infer that $I$ is a positive
linear functional on $\mathcal{C}(X)$, so, in view of Riesz representation
theorem, we conclude that there exists a unique borelian positive measure $%
\mu _{\mathcal{S}}$ on\textit{\ }$X$ such that $c_{g}=\underset{X}{\int }%
gd\mu _{\mathcal{S}}$ for every $g\in \mathcal{C}(X)$. $\square $

\bigskip

\textbf{Theorem 4.9.\ }$M_{\mathcal{S}}:(\mathcal{M}(X),d_{H})\rightarrow (%
\mathcal{M}(X),d_{H})$ \textit{is a Picard operator for every} $\varphi $-$%
\max $-\textit{IFSp} $\mathcal{S}=((X,d),(f_{i})_{i\in
\{1,2,...,m\}},(p_{i})_{i\in \{1,2,...,m\}})$, \textit{with }$(X,d)$\textit{%
\ compact}, \textit{and the support of the fixed point of} $M_{\mathcal{S}}$ 
\textit{is} $A_{\mathcal{S}}$\textit{.}

\textit{Proof}. First of all let us note that Remark 3.6 can be restated as $%
\underset{X}{\int }gdM_{\mathcal{S}}(\nu )$\linebreak $=\underset{X}{\int }%
B_{\mathcal{S}}(g)d\nu $ for every $g\in \mathcal{C}(X)$ and every $\nu \in 
\mathcal{M}(X)$. Therefore, we get%
\begin{equation}
\underset{X}{\int }gdM_{\mathcal{S}}^{[n]}(\nu )=\underset{X}{\int }B_{%
\mathcal{S}}^{[n]}(g)d\nu \text{,}  \tag{1}
\end{equation}%
for every $g\in \mathcal{C}(X)$, $\nu \in \mathcal{M}(X)$ and $n\in \mathbb{N%
}$.

Now we divide the proof into three steps.

\medskip

\textbf{Step 1}. \textit{The measure} $\mu _{\mathcal{S}}$,\textit{\
provided by Proposition 4.8, belongs to }$\mathcal{M}(X)$.

For the function $g_{0}:X\rightarrow \mathbb{R}$ given by $g_{0}(x)=1$ for
every $x\in X$, we have $B_{\mathcal{S}}^{[n]}g_{0}(x)=1$ for every $x\in X$
and every $n\in \mathbb{N}$. Consequently we have $c_{g_{0}}\overset{\text{%
Proposition 4.7}}{=}1$, i.e. $\underset{X}{\int }1d\mu _{\mathcal{S}}=1$, so 
$\mu _{\mathcal{S}}(X)=1$. Moreover, as supp $\mu _{\mathcal{S}}$ is a
closed subset of the compact set $X$, it is compact.

\medskip

\textbf{Step 2}. \textit{The measure} $\mu _{\mathcal{S}}\in \mathcal{M}(X)$,%
\textit{\ provided by Proposition 4.8, is the unique fixed point of} $M_{%
\mathcal{S}}$ \textit{and} supp $\mu _{\mathcal{S}}=A_{S}$\textit{.}

On the one hand, by passing to limit as $n\rightarrow \infty $ in the
relation $B_{\mathcal{S}}^{[n+1]}(g)=B_{\mathcal{S}}^{[n]}(B_{\mathcal{S}%
}(g))$ which is valid for every $n\in \mathbb{N}$ and every $g\in \mathcal{C}%
(X)$, taking into account Proposition 4.7, we get $\underset{X}{\int }gd\mu
_{\mathcal{S}}=\underset{X}{\int }B_{\mathcal{S}}(g)d\mu _{\mathcal{S}}$ for
every $g\in \mathcal{C}(X)$. In other words, we have $\underset{X}{\int }%
gd(\mu _{\mathcal{S}}-M_{\mathcal{S}}(\mu _{\mathcal{S}}))=0$ for every $%
g\in \mathcal{C}(X)$ which implies that $M_{\mathcal{S}}(\mu _{\mathcal{S}%
})=\mu _{\mathcal{S}}$, i.e. $\mu _{\mathcal{S}}$ is a fixed point of $M_{%
\mathcal{S}}$. Since supp $\mu _{\mathcal{S}}=$ supp $M_{\mathcal{S}}(\mu _{%
\mathcal{S}})\overset{\text{Lemma 3.4}}{=}F_{\mathcal{S}}($supp $\mu _{%
\mathcal{S}})$, we infer that supp $\mu _{\mathcal{S}}$ is the fixed point
of $F_{\mathcal{S}}$, so supp $\mu _{\mathcal{S}}=A_{S}$\textit{.}

On the other hand, if $\nu \in \mathcal{M}(X)$ has the property that $M_{%
\mathcal{S}}(\nu )=\nu $, then, from $(1)$, we get $\underset{X}{\int }gd\nu
=\underset{X}{\int }B_{\mathcal{S}}^{[n]}(g)d\nu $ for every $g\in \mathcal{C%
}(X)$ and every $n\in \mathbb{N}$. By passing to limit as $n\rightarrow
\infty $, based on Proposition 4.7, we get $\underset{X}{\int }gd\nu =%
\underset{X}{\int }c_{g}d\nu $. As $\nu (X)=1$, we obtain $\underset{X}{\int 
}gd\nu =\underset{X}{\int }gd\mu _{\mathcal{S}}$, so $\underset{X}{\int }%
gd(\nu -\mu _{\mathcal{S}})=0$ for every $g\in \mathcal{C}(X)$. We conclude
that $\nu =\mu _{\mathcal{S}}$, i.e. $\mu _{\mathcal{S}}$ is the unique
fixed point of $M_{\mathcal{S}}$.

\medskip

\textbf{Step 3}. $\underset{n\rightarrow \infty }{\lim }M_{\mathcal{S}%
}^{[n]}(\nu )=\mu _{\mathcal{S}}$ \textit{for every }$\nu \in \mathcal{M}(X)$%
.

Since, based on Proposition 4.7, we have $\underset{n\rightarrow \infty }{%
\lim }\underset{X}{\int }B_{\mathcal{S}}^{[n]}(g)d\nu =\underset{X}{\int }%
c_{g}d(\nu )$, as $\nu (X)=1$, we infer that $\underset{n\rightarrow \infty }%
{\lim }\underset{X}{\int }B_{\mathcal{S}}^{[n]}(g)d\nu =\underset{X}{\int }%
gd\mu _{\mathcal{S}}$. Hence, using $(1)$, we get $\underset{n\rightarrow
\infty }{\lim }\underset{X}{\int }gdM_{\mathcal{S}}^{[n]}(\nu )=\underset{X}{%
\int }gd\mu _{\mathcal{S}}$ for every $g\in \mathcal{C}(X)$. Taking into
account Remark 2.2, we conclude that $\underset{n\rightarrow \infty }{\lim }%
d_{H}(M_{\mathcal{S}}^{[n]}(\nu ),\mu _{\mathcal{S}})=0$.

\medskip

The last two steps assure us that $M_{\mathcal{S}}$ is a Picard operator. $%
\square $

\bigskip

\textbf{B. The case of a general} $\varphi $-$\max $\textbf{-IFSp}

\bigskip

Let us start with some:

\medskip

\textbf{Notations}.

For a complete metric space$\mathcal{\ }(X,d)$ and a compact subset $Y$ of $%
X $, we shall consider the following:

$\bullet $ the function $R_{Y}:Lip_{1}(X,\mathbb{R})\rightarrow Lip_{1}(Y,%
\mathbb{R})$ given by $R_{Y}(f)(x)=f(x)$ for every $f\in Lip_{1}(X,\mathbb{R}%
)$ and every $x\in Y$

$\bullet $ the function $E_{Y}:Lip_{1}(T,\mathbb{R})\rightarrow Lip_{1}(X,%
\mathbb{R})$ given by $E_{Y}(f)(x)=\underset{y\in Y}{\sup }%
(f(y)-lip(f)d(x,y))$ for every $f\in Lip_{1}(Y,\mathbb{R})$ and every $x\in
X $

$\bullet $ the function $i_{Y}:\mathcal{M}(Y)\rightarrow \mathcal{M}(X)$
given by $i_{Y}(\mu )(B)=\mu (Y\cap B)$ for every $\mu \in \mathcal{M}(Y)$
and every $B\in \mathcal{B}(X)$

$\bullet $ the Hutchinson distance $d_{H}^{Y}:\mathcal{M}(Y)\times \mathcal{M%
}(Y)\rightarrow \lbrack 0,\infty )$ described by $d_{H}^{Y}(\mu ,\nu )=%
\underset{f\in Lip_{1}(Y,\mathbb{R})}{\sup }\left\vert \underset{Y}{\dint }%
fd\mu -\underset{Y}{\dint }fd\nu \right\vert $\textit{\ }for every $\mu ,\nu
\in \mathcal{M}(Y)$.

For a $\varphi $-$\max $-IFSp $\mathcal{S}=((X,d),(f_{i})_{i\in
\{1,2,...,m\}},(p_{i})_{i\in \{1,2,...,m\}})$ and $\nu \in \mathcal{M}(X)$
we shall consider the set $K_{\nu }\overset{def}{=}A_{\mathcal{S}}\cup (%
\underset{n\in \mathbb{N}}{\cup }F_{\mathcal{S}}^{[n]}($supp $\nu ))$.

\bigskip

Note that the functions $R_{Y}$ and $E_{Y}$ are well-defined (i.e. $%
R_{Y}(f)\in Lip_{1}(Y,\mathbb{R})$ for every $f\in Lip_{1}(X,\mathbb{R})$
and $E_{Y}(f)\in Lip_{1}(X,\mathbb{R})$ for every $f\in Lip_{1}(Y,\mathbb{R}%
) $; moreover, $R_{Y}(E_{Y}(f))=f$ and $lip(E_{Y}(f))=lip(f)$, according to
a famous result due to E.J. McShane -see Theorem 1 from [13]-).

\bigskip

\textbf{Remark 4.10.} \textit{Given a complete metric space }$(X,d)$\textit{%
\ and a compact subset }$Y$\textit{\ of }$X$\textit{, for every} $\mu \in 
\mathcal{M}(Y)$, \textit{we have:}

\textit{i) }$\underset{Y}{\int }R_{Y}(f)d\mu =\underset{X}{\int }%
fd(i_{Y}(\mu ))$\textit{\ for every }$f\in Lip_{1}(X,\mathbb{R})$\textit{;}

\textit{ii) }$\underset{Y}{\int }fd\mu =\underset{X}{\int }%
E_{Y}(f)d(i_{Y}(\mu ))$\textit{\ for every }$f\in Lip_{1}(Y,\mathbb{R})$%
\textit{.}

Indeed, since supp $i_{Y}(\mu )=\underset{F=\overline{F}\subseteq X\text{, }%
i_{Y}(\mu )(F)=i_{Y}(\mu )(X)}{\cap }F=\underset{F=\overline{F}\subseteq X%
\text{, }\mu (F\cap Y)=\mu (Y)}{\cap }F$\linebreak $\subseteq Y$ (as $Y=%
\overline{Y}\subseteq X$ and $\mu (Y\cap Y)=\mu (Y)$), we have $\underset{X}{%
\int }fd(i_{Y}(\mu ))=\underset{Y}{\int }fd(i_{Y}(\mu ))=\underset{Y}{\int }%
R_{Y}(f)d\mu $ for every $f\in Lip_{1}(X,\mathbb{R})$ and $\underset{X}{\int 
}E_{Y}(f)d(i_{Y}(\mu ))=\underset{Y}{\int }E_{Y}(f)d(i_{Y}(\mu ))=\underset{Y%
}{\int }fd\mu $ for every $f\in Lip_{1}(Y,\mathbb{R})$.

\bigskip

\textbf{Lemma 4.11.} \textit{Given a complete metric space }$(X,d)$\textit{\
and a compact subset }$Y$\textit{\ of }$X$\textit{,} \textit{we have }$%
d_{H}^{Y}(\mu _{1},\mu _{2})=d_{H}(i_{Y}(\mu _{1}),i_{Y}(\mu _{2}))$\textit{%
\ for every} $\mu _{1},\mu _{2}\in \mathcal{M}(Y)$.

\textit{Proof}. On the one hand, we have $d_{H}^{Y}(\mu _{1},\mu _{2})=%
\underset{f\in Lip_{1}(Y,\mathbb{R})}{\sup }\left\vert \underset{Y}{\dint }%
fd\mu _{1}-\underset{Y}{\dint }fd\mu _{2}\right\vert $\linebreak $\overset{%
\text{Remark 4.10, ii)}}{=}\underset{f\in Lip_{1}(Y,\mathbb{R})}{\sup }%
\left\vert \underset{X}{\dint }E_{Y}(f)d(i_{Y}(\mu _{1}))-\underset{X}{\dint 
}E_{Y}(f)d(i_{Y}(\mu _{2}))\right\vert \leq $\linebreak $\underset{f\in
Lip_{1}(X,\mathbb{R})}{\sup }\left\vert \underset{X}{\dint }fd(i_{Y}(\mu
_{1}))-\underset{X}{\dint }fd(i_{Y}(\mu _{2}))\right\vert =d_{H}(i_{Y}(\mu
_{1}),i_{Y}(\mu _{2}))$, so%
\begin{equation}
d_{H}^{Y}(\mu _{1},\mu _{2})\leq d_{H}(i_{Y}(\mu _{1}),i_{Y}(\mu _{2}))\text{%
,}  \tag{1}
\end{equation}%
for every $\mu _{1},\mu _{2}\in \mathcal{M}(Y)$.

On the other hand $d_{H}(i_{Y}(\mu _{1}),i_{Y}(\mu _{2}))=\underset{f\in
Lip_{1}(X,\mathbb{R})}{\sup }\left\vert \underset{X}{\dint }fd(i_{Y}(\mu
_{1}))-\underset{X}{\dint }fd(i_{Y}(\mu _{2}))\right\vert $\linebreak $%
\overset{\text{Remark 4.10, i)}}{=}\underset{f\in Lip_{1}(X,\mathbb{R})}{%
\sup }\left\vert \underset{X}{\dint }R_{Y}(f)d\mu _{1}-\underset{X}{\dint }%
R_{Y}(f)d\mu _{2}\right\vert \leq \underset{f\in Lip_{1}(Y,\mathbb{R})}{\sup 
}\left\vert \underset{Y}{\dint }fd\mu _{1}-\underset{Y}{\dint }fd\mu
_{2}\right\vert $\linebreak $=d_{H}^{Y}(\mu _{1},\mu _{2})$, so%
\begin{equation}
d_{H}(i_{Y}(\mu _{1}),i_{Y}(\mu _{2}))\leq d_{H}^{Y}(\mu _{1},\mu _{2})\text{%
,}  \tag{2}
\end{equation}%
for every $\mu _{1},\mu _{2}\in \mathcal{M}(Y)$.

From $(1)$ and $(2)$, we get the conclusion. $\square $

\bigskip

\textbf{Remark 4.12.} \textit{For every }$\varphi $\textit{-}$\max $\textit{%
-IFSp} $\mathcal{S}=((X,d),(f_{i})_{i\in \{1,...,m\}},(p_{i})_{i\in
\{1,...,m\}})$ \textit{and} $\nu \in \mathcal{M}(X)$, \textit{we have:}

\textit{i)} $F_{\mathcal{S}}(K_{\nu })\subseteq K_{\nu }$\textit{;}

\textit{ii) }$K_{\nu }$\textit{\ is compact.}

Indeed, we have $F_{\mathcal{S}}(K_{\nu })=F_{\mathcal{S}}(A_{\mathcal{S}%
}\cup (\underset{n\in \mathbb{N}}{\cup }F_{\mathcal{S}}^{[n]}($supp $\nu
)))\subseteq F_{\mathcal{S}}(A_{\mathcal{S}})\cup (\underset{n\in \mathbb{N}}%
{\cup }F_{\mathcal{S}}(F_{\mathcal{S}}^{[n]}($supp $\nu ))))=F_{\mathcal{S}%
}(A_{\mathcal{S}})\cup (\underset{n\in \mathbb{N}}{\cup }(F_{\mathcal{S}%
}^{[n+1]}($supp $\nu )))\subseteq K_{\nu }$. The compactness of $K_{\nu }$
is assured, via Proposition 2.5, by the fact that\linebreak\ $\underset{%
n\rightarrow \infty }{\lim }H(F_{\mathcal{S}}^{[n]}($supp $\nu )),A_{%
\mathcal{S}})=0$.

\bigskip

Taking into account that $f_{i}(A_{\mathcal{S}})\subseteq A_{\mathcal{S}}$
and $f_{i}(K_{\nu })\subseteq K_{\nu }$ for every $i\in \{1,2,...,m\}$ and
every $\nu \in \mathcal{M}(X)$, we can consider the $\varphi $\textit{-}$%
\max $-IFSps $\mathcal{S}_{A_{\mathcal{S}}}=((A_{\mathcal{S}},d),(\phi
_{i})_{i\in \{1,...,m\}},(p_{i})_{i\in \{1,...,m\}})$ and $\mathcal{S}%
_{K_{\nu }}=((K_{\nu },d),(\psi _{i})_{i\in \{1,...,m\}},(p_{i})_{i\in
\{1,...,m\}})$, where $\phi _{i}(x)=f_{i}(x)$ for every $x\in A_{\mathcal{S}%
} $ and every $i\in \{1,2,...,m\}$ and $\psi _{i}(x)=f_{i}(x)$ for every $%
x\in K_{\nu }$ and every $i\in \{1,2,...,m\}$. We can also consider the
Markov operator $M_{\mathcal{S}}^{A_{\mathcal{S}}}:\mathcal{M}(A_{\mathcal{S}%
})\rightarrow \mathcal{M}(A_{\mathcal{S}})$ associated to\ $\mathcal{S}_{A_{%
\mathcal{S}}}$ and the Markov operator $M_{\mathcal{S}}^{K_{\nu }}:\mathcal{M%
}(K_{\nu })\rightarrow \mathcal{M}(K_{\nu })$ associated to\ $\mathcal{S}%
_{_{K_{\nu }}}$. According to Theorem 4.9, $M_{\mathcal{S}}^{A_{\mathcal{S}%
}} $ and $M_{\mathcal{S}}^{K_{\nu }}$ are Picard operators and we denote the
fixed point of $M_{\mathcal{S}}^{A_{\mathcal{S}}}$ by $\mu _{\mathcal{S}%
}^{A_{\mathcal{S}}}$ and the fixed point of $M_{\mathcal{S}}^{K_{\nu }}$ by $%
\mu _{\mathcal{S}}^{K_{\nu }}$.

\bigskip

\textbf{Proposition 4.13.} \textit{The Markov operator }$M_{\mathcal{S}}$ 
\textit{associated to a }$\varphi $\textit{-}$\max $\textit{-IFSp} $\mathcal{%
S}$ \textit{has a unique fixed point denoted by} $\mu _{\mathcal{S}}$ 
\textit{whose support is} $A_{\mathcal{S}}$.

\textit{Proof}. The function $G_{\mathcal{S}}:\{\mu \in \mathcal{M}(A_{%
\mathcal{S}})\mid M_{\mathcal{S}}^{A_{\mathcal{S}}}(\mu )=\mu \}\rightarrow
\{\mu \in \mathcal{M}(X)\mid M_{\mathcal{S}}(\mu )=\mu \}$ given by $G_{%
\mathcal{S}}(\mu )=i_{A_{\mathcal{S}}}(\mu )$ for every $\mu \in \mathcal{M}%
(A_{\mathcal{S}})$ such that $M_{\mathcal{S}}^{A_{\mathcal{S}}}(\mu )=\mu $
is a bijection (whose inverse is the function $H_{\mathcal{S}}:\{\mu \in 
\mathcal{M}(X)\mid M_{\mathcal{S}}(\mu )=\mu \}\rightarrow \{\mu \in 
\mathcal{M}(A_{\mathcal{S}})\mid M_{\mathcal{S}}^{A_{\mathcal{S}}}(\mu )=\mu
\}$ given by $H_{\mathcal{S}}(\mu )=\mu _{\mid \mathcal{B}(A_{\mathcal{S}})}$
for every $\mu \in \mathcal{M}(X)$ such that $M_{\mathcal{S}}(\mu )=\mu $.
Therefore $i_{A_{\mathcal{S}}}(\mu _{\mathcal{S}}^{A_{\mathcal{S}}})\overset{%
not}{=}\mu _{\mathcal{S}}$ is the unique fixed point of $M_{\mathcal{S}}$.
In addition, supp $\mu _{\mathcal{S}}=$supp $M_{\mathcal{S}}(\mu _{\mathcal{S%
}})\overset{\text{Lemma 3.4}}{=}F_{\mathcal{S}}($supp $\mu _{\mathcal{S}})$,
so, taking into account the uniqueness of the fixed point of $F_{\mathcal{S}%
} $, we infer that supp $\mu _{\mathcal{S}}=A_{\mathcal{S}}$. $\square $

\bigskip

\textbf{Lemma 4.14.} \textit{For every }$\varphi $\textit{-}$\max $\textit{%
-IFSp} $\mathcal{S}=((X,d),(f_{i})_{i\in \{1,...,m\}},(p_{i})_{i\in
\{1,...,m\}})$ \textit{and} $\nu \in \mathcal{M}(X)$\textit{, the measure} $%
\nu _{0}$\textit{, given by} $\nu _{0}(B)=\nu (B\cap K_{\nu })$ \textit{for
every} $B\in \mathcal{B}(K_{\nu })=\{B^{^{\prime }}\cap K_{\nu }\mid
B^{^{\prime }}\in \mathcal{B}(X)\}$\textit{, has the following properties:}

\textit{i)} $\nu _{0}\in \mathcal{M}(K_{\nu })$\textit{;}

\textit{ii) }$i_{K_{\nu }}(\nu _{0})=\nu $\textit{.}

\textit{Proof}. We start by noting that, as $B^{^{\prime }}\smallsetminus
K_{\nu }\subseteq B^{^{\prime }}\smallsetminus $supp $\nu \subseteq
X\smallsetminus $supp $\nu $, we have $0\leq \nu (B^{^{\prime
}}\smallsetminus K_{\nu })\leq \nu (B^{^{\prime }}\smallsetminus $supp $\nu
)\leq \nu (X\smallsetminus $supp $\nu )\overset{\text{definition of supp}}{=}%
0$, so%
\begin{equation}
\nu (B^{^{\prime }}\smallsetminus K_{\nu })=0\text{,}  \tag{1}
\end{equation}%
for every $B^{^{\prime }}\in \mathcal{B}(X)$.

Then $1\overset{\nu \in \mathcal{M}(X)}{=}\nu (X)=\nu (X\smallsetminus
K_{\nu })+\nu (K_{\nu })\overset{(1)}{=}\nu (K_{\nu })$, so i) is proved.

Moreover, $i_{K_{\nu }}(\nu _{0})(B^{^{\prime }})=\nu _{0}(B^{^{\prime
}}\cap K_{\nu })=\nu (B^{^{\prime }}\cap K_{\nu })\overset{(1)}{=}\nu
(B^{^{\prime }}\cap K_{\nu })+$\linebreak $\nu (B^{^{\prime }}\smallsetminus
K_{\nu })=\nu (B^{^{\prime }})$ for every $B^{^{\prime }}\in \mathcal{B}(X)$%
, so ii) is also proved. $\square $

\bigskip

\textbf{Lemma 4.15.} \textit{Given a }$\varphi $\textit{-}$\max $\textit{%
-IFSp} $\mathcal{S}=((X,d),(f_{i})_{i\in \{1,...,m\}},(p_{i})_{i\in
\{1,...,m\}})$\textit{, we have }$i_{K_{\nu }}((M_{\mathcal{S}}^{K_{\nu
}})^{[n]}(\nu ))=M_{\mathcal{S}}^{[n]}(i_{K_{\nu }}(\nu ))$ \textit{for every%
} $n\in \mathbb{N}$ \textit{and every} $\nu \in \mathcal{M}(K_{\nu })$.

\textit{Proof}. First we prove that%
\begin{equation}
i_{K_{\nu }}(M_{\mathcal{S}}^{K_{\nu }}(\nu ))=M_{\mathcal{S}}(i_{K_{\nu
}}(\nu ))\text{,}  \tag{1}
\end{equation}%
for every $\nu \in \mathcal{M}(K_{\nu })$.

In order to justify $(1)$, it suffices to check that $\underset{X}{\dint }%
gd(i_{K_{\nu }}(M_{\mathcal{S}}^{K_{\nu }}(\nu )))=p_{1}\underset{X}{\dint }%
g\circ f_{1}d(i_{K_{\nu }}(\nu ))+...+p_{m}\underset{X}{\dint }g\circ
f_{m}d(i_{K_{\nu }}(\nu ))$ for every continuous and bounded function $%
g:X\rightarrow \mathbb{R}$. This is true since $\underset{X}{\dint }%
gd(i_{K_{\nu }}(M_{\mathcal{S}}^{K_{\nu }}(\nu )))\overset{\text{Remark
4.10, i)}}{=}\underset{K_{\nu }}{\dint }R_{K_{\nu }}(g)d(M_{\mathcal{S}%
}^{K_{\nu }}(\nu ))=p_{1}\underset{K_{\nu }}{\dint }R_{K_{\nu }}(g)\circ
\psi _{1}d\nu +...+p_{m}\underset{K_{\nu }}{\dint }R_{K_{\nu }}(g)\circ \psi
_{m}d\nu =p_{1}\underset{K_{\nu }}{\dint }R_{K_{\nu }}(g\circ f_{1})d\nu
+...+p_{m}\underset{K_{\nu }}{\dint }R_{K_{\nu }}(g\circ f_{m})d\nu \overset{%
\text{Remark 4.10, ii)}}{=}p_{1}\underset{X}{\dint }g\circ f_{1}d(i_{K_{\nu
}}(\nu ))+...+p_{m}\underset{X}{\dint }g\circ f_{m}d(i_{K_{\nu }}(\nu ))$.

Now, for every $n\in \mathbb{N}$, we have $i_{K_{\nu }}((M_{\mathcal{S}%
}^{K_{\nu }})^{[n]}(\nu ))=i_{K_{\nu }}((M_{\mathcal{S}}^{K_{\nu }}(M_{%
\mathcal{S}}^{K_{\nu }})^{[n-1]}(\nu )))$\linebreak $\overset{\text{(1)}}{=}%
M_{\mathcal{S}}(i_{K_{\nu }}((M_{\mathcal{S}}^{K_{\nu }})^{[n-1]}(\nu )))=M_{%
\mathcal{S}}^{[2]}(i_{K_{\nu }}((M_{\mathcal{S}}^{K_{\nu }})^{[n-2]}(\nu
)))=...=M_{\mathcal{S}}^{[n]}(i_{K_{\nu }}(\nu ))$ for every $\nu \in 
\mathcal{M}(K_{\nu })$. $\square $

\bigskip

\textbf{Lemma 4.16.} \textit{Given a }$\varphi $\textit{-}$\max $\textit{%
-IFSp} $\mathcal{S}=((X,d),(f_{i})_{i\in \{1,...,m\}},(p_{i})_{i\in
\{1,...,m\}})$\textit{, we have }$\mu _{\mathcal{S}}=i_{K_{\nu }}(\mu _{%
\mathcal{S}}^{K_{\nu }})$\textit{\ for every} $\nu \in \mathcal{M}(X)$.

\textit{Proof}. As in the proof of Proposition 4.13, one can check that the
function $\mu \mapsto i_{K_{\nu }}(\mu )$ is a bijection from the fixed
points of $M_{\mathcal{S}}^{K_{\nu }}$ to the fixed point of $M_{\mathcal{S}%
} $, so $\mu _{\mathcal{S}}=i_{K_{\nu }}(\mu _{\mathcal{S}}^{K_{\nu }})$. $%
\square $

\bigskip

\textbf{Proposition 4.17.} \textit{Given a }$\varphi $\textit{-}$\max $%
\textit{-IFSp} $\mathcal{S}=((X,d),(f_{i})_{i\in \{1,...,m\}},(p_{i})_{i\in
\{1,...,m\}})$\textit{, we have }$\underset{n\rightarrow \infty }{\lim }%
d_{H}(M_{\mathcal{S}}^{[n]}(\nu ),\mu _{\mathcal{S}})=0$\textit{\ for every} 
$\nu \in \mathcal{M}(X)$.

\textit{Proof}. According to Lemma 4.14, for every $\nu \in \mathcal{M}(X)$
there exists $\nu _{0}\in \mathcal{M}(K_{\nu })$\textit{\ }such that\textit{%
\ }$i_{K_{\nu }}(\nu _{0})=\nu $. Then we have $d_{H}(M_{\mathcal{S}%
}^{[n]}(\nu ),\mu _{\mathcal{S}})=d_{H}(M_{\mathcal{S}}^{[n]}(i_{K_{\nu
}}(\nu _{0})),\mu _{\mathcal{S}})\overset{\text{Lemma 4.15}}{=}%
d_{H}(i_{K_{\nu }}((M_{\mathcal{S}}^{K_{\nu }})^{[n]}(\nu _{0})),\mu _{%
\mathcal{S}})\overset{\text{Lemma 4.16}}{=}$\linebreak $d_{H}(i_{K_{\nu
}}((M_{\mathcal{S}}^{K_{\nu }})^{[n]}(\nu _{0})),i_{K_{\nu }}(\mu _{\mathcal{%
S}}^{K_{\nu }}))\overset{\text{Lemma 4.11}}{=}d_{H}^{K_{\nu }}((M_{\mathcal{S%
}}^{K_{\nu }})^{[n]}(\nu _{0})),\mu _{\mathcal{S}}^{K_{\nu }})$ for every $%
n\in \mathbb{N}$, so $\underset{n\rightarrow \infty }{\lim }d_{H}(M_{%
\mathcal{S}}^{[n]}(\nu ),\mu _{\mathcal{S}})=\underset{n\rightarrow \infty }{%
\lim }d_{H}^{K_{\nu }}((M_{\mathcal{S}}^{K_{\nu }})^{[n]}(\nu _{0})),\mu _{%
\mathcal{S}}^{K_{\nu }})\overset{\text{Theorem 4.9}}{=}0$. $\square $

\bigskip

Combining Proposition 4.13 with Proposition 4.17, we get the following

\bigskip

\textbf{Theorem 4.18.\ }$M_{\mathcal{S}}:(\mathcal{M}(X),d_{H})\rightarrow (%
\mathcal{M}(X),d_{H})$ \textit{is a Picard operator for every} $\varphi $-$%
\max $-\textit{IFSp} $\mathcal{S}=((X,d),(f_{i})_{i\in
\{1,2,...,m\}},(p_{i})_{i\in \{1,2,...,m\}})$ \textit{and the support of the
fixed point of} $M_{\mathcal{S}}$ \textit{is} $A_{\mathcal{S}}$\textit{.}

\bigskip

\textbf{References}

\bigskip

[1] Barnsley, M.F., Demko, S.G., Elton, J.H. and Geromino, J.S.: Invariant
measures for Markov processes arising from iterated function systems with
place-dependent probabilities, Ann. Inst. H. Poincar\'{e} Probab. Statist., 
\textbf{24 }(1988), 367--394 and \textbf{25 }(1989), 589--590.

[2] Barnsley, M.F and Elton, J.H.: A new class of Markov processes for image
encoding, Adv. in Appl. Prob., \textbf{20 }(1988), 14--32.

[3] Chi\c{t}escu I. and Ni\c{t}\u{a} L.: Fractal vector measures, Sci.
Bull., Ser. A, Appl. Math. Phys., Politeh. Univ. Buchar., \textbf{77}
(2015), 219-228.

[4] Elton, J.H.: An ergodic theorem for iterated maps, Ergod. Theory Dynam.
Systems, \textbf{7 }(1987), 481--488.

[5] Fisher, Y.: Fractal Image Compression; Theory and Application,
Springer-Verlag, New York, 1995.

[6] Georgescu, F., Miculescu, R. and Mihail, A.: Iterated function systems
consisting of $\varphi $-max-contractions have attractor, available at
arXiv:1704.02652.

[7] Hutchinson, J.E.: Fractals and self similarity, Indiana Univ. Math. J., 
\textbf{30} (1981), 713--747.

[8] Iosifescu, M. and Grigorescu, S.: Dependence with complete connections
and its applications, Cambridge Univ. Press, Cambridge, 1990.

[9] Jaroszewska, J.: Iterated function systems with continuous place
dependent probabilities, Univ. Iagell. Acta Math., \textbf{40} (2002),
137-146.

[10] Jaroszewska, J.: How to construct asymptotically stable iterated
function systems, Stat. Probab. Lett., \textbf{78 }(2008), 1570-1576.

[11] Lasota, A.\ and Myjak, J.: Markov operators and fractals, Bull. Pol.
Acad. Sci., Math., \textbf{45 }(1997), 197-210.

[12] Marchi, M.V.: Invariant measures in quasi-metric spaces, Z. Anal.
Anwend., \textbf{22} (2003), 17--32.

[13] McShane, E.J.: Extension of range of functions, Bull. Amer. Math. Soc., 
\textbf{40} (1934), 837-842.

[14] Mendivil, F.: A generalization of IFS with probabilities to infinitely
many maps, Rocky Mt. J. Math., \textbf{28} (1998), 1043--1051.

[15] Miculescu, R.: Generalized iterated function systems with place
dependent probabilities, Acta. Appl. Math., \textbf{130} (2014), 135-150.

[16] Miculescu, R. and Mihail, A.: A generalization of Istr\u{a}\c{t}escu's
fixed point theorem for convex contractions, Fixed Point Theory, in print,
arXiv:1512.05490.

[17] Mihail, A. and Miculescu, R.: A generalization of the Hutchinson
measure. Mediterr. J. Math., \textbf{6}\ (2009), 203--213.

[18] Mosco, U.: Self-similar measures in quasi-metric spaces. In:
Festschrift Dedicated to Alfonso Vignoli on the Occasion of His 60th
Birthday. Prog. Nonlinear Differ. Equ. Appl., vol. 40, pp. 233--248. Birkh%
\"{a}user, Basel, 2000.

[19] \"{O}berg, A.: Approximation of invariant measures for random
iterations, Rocky Mt. J. Math., \textbf{36 }(2006), 273-301.

[20] Oliveira, O.: The ergodic theorem for a new kind of attractor of a
GIFS, Chaos Solitons Fractals, \textbf{98} (2017), 63-71.

[21] Onicescu, O. and Mihoc, G.: Sur les cha\^{\i}nes de variables
statisques, Bull. Soc. Math. Fr., \textbf{59} (1935), 174--192.

[22] Secelean, N.: Countable Iterated Function Systems, LAP Lambert Academic
Publishing, 2013.

[23] Secelean, N.: Invariant measure associated with a generalized countable
iterated function system, Mediterr. J. Math.. \textbf{11} (2014), 361-372.

[24] Stenflo \"{O}.: Ergodic theorems for Markov chains represented by
iterated function systems, Bull. Pol. Acad. Sci., Math., \textbf{49} (2001),
27-43.

[25] Stenflo \"{O}.: Uniqueness of Invariant Measures for Place-Dependent
Random Iterations of Functions. In: Barnsley M.F., Saupe D., Vrscay E.R.
(eds) Fractals in Multimedia. The IMA Volumes in Mathematics and its
Application, vol 132. Springer, New York, NY, 2002.

[26] Szarek, T.: Lower bound technique and its applications to function
systems and stochastic partial differential equations, Math. Appl. (Warsaw), 
\textbf{41 }(2013), 185-198.

[27] Werner, I.: Contractive Markov systems, J. Lond. Math. Soc., \textbf{71}
(2005), 236-258.

[28] Zaharopol, R.: Equicontinuity and existence of attractive probability
measures for some iterated function systems, Rev. Roum. Math. Pures Appl., 
\textbf{52} (2007), 259-286.

\bigskip

{\small Flavian Georgescu}

{\small Faculty of Mathematics and Computer Science}

{\small University of Pite\c{s}ti, Romania}

{\small T\^{a}rgul din Vale 1, 110040, Pite\c{s}ti, Arge\c{s}, Romania}

{\small E-mail: faviu@yahoo.com}

\bigskip

{\small Radu Miculescu}

{\small Faculty of Mathematics and Computer Science}

{\small Bucharest University, Romania}

{\small Str. Academiei 14, 010014, Bucharest}

{\small E-mail: miculesc@yahoo.com}

\bigskip

{\small Alexandru Mihail}

{\small Faculty of Mathematics and Computer Science}

{\small Bucharest University, Romania}

{\small Str. Academiei 14, 010014, Bucharest}

{\small E-mail: mihail\_alex@yahoo.com}

\end{document}